\documentclass{notices}
\usepackage{amsfonts,amssymb,amsmath,amscd}

\usepackage[all]{xy}
\usepackage{url}
\usepackage{stmaryrd}
\usepackage{color}
\usepackage{graphicx}
\usepackage{caption}

\newtheorem{thm}{Theorem}

\newtheorem{defn}[thm]{Definition}

\newtheorem{ex}[thm]{Example}

\numberwithin{thm}{section}

\newcommand{\ga}[2]{\begin{gather}\label{#1}#2 \end{gather}}

\newcommand{\sO}{{\mathcal O}}

\newcommand{\A}{{\mathbb A}}

\newcommand{\C}{{\mathbb C}}

\newcommand{\F}{{\mathbb F}}

\renewcommand{\P}{{\mathbb P}}
\newcommand{\Q}{{\mathbb Q}}

\newcommand{\Z}{{\mathbb Z}}

\DeclareMathOperator{\GL}{GL}

\title{
Some Arithmetic Properties of Complex Local Systems
}

\author{
  H\'el\`ene Esnault
  \affil{Freie Universit\"at Berlin, 
  Mathematik\\
  Arnimallee 7,
  14195 Berlin, 
  Germany
  }
  \affil{esnault@math.fu-berlin.de
   }
  }

\begin{document}

\maketitle

\section{Introduction} \label{sec:intro}
A group  $\pi$ is said to be {\it finitely generated} if it is spanned by finitely many letters, that is, if it is the quotient $F\to \pi$ of a free group $F$ on finitely many letters. It is said to be {\it finitely presented} if the kernel of such a quotient is itself finitely generated. This does not depend on the choice of generation chosen.  For example the trivial group $\pi=\{1\}$ is surely finitely presented as the quotient of the free group in $1$ generator by itself (!). The following finitely presented group shall play a role in the note:
\begin{ex} \label{ex:1}
 The group  $\Gamma_0$ is generated  by two elements $(a,b)$ with one relation $b^2=a^2ba^{-2}$.
  \end{ex}
  There are groups which are finitely generated but not finitely presented, see the interesting MathOverflow elementary discussion on the topic (\url{https://tinyurl.com/3cavr69a}).

\medskip

The finitely presented groups appear naturally in many branches of mathematics. The fundamental group $\pi_1(M,m)$ of a topological space $M$ based at a point $m$ is defined to be the group of homotopy classes of loops centered at $m$. 
A group is finitely presented {\it if and only if}  
  it is the fundamental group  $\pi_1(M,m)$ {\it of  a connected  finite $CW$-complex $M$ }based at a point $m$.     This is essentially by definition of a $CW$ ($C=$closure-finite,  $W=$weak)  complex which is a
 topological space defined by an increasing sequence of topological subspaces, each one obtained by  gluing cells of growing dimension to the previous one.  So the $1$-cells glued to the $0$-cell  $m$ yield  the loops  on which we take the free group $F$, and the relations come from the finitely many $2$-cells glued to the loops.

\medskip

   If $X$ is a smooth connected  quasi-projective  complex 
variety,  its complex points $X(\mathbb C)$ form a topological manifold which   has the homotopy type of  a connected  finite $CW$-complex $M$. 

\medskip

 The difference between $X$ and its complex points $X(\C)$  is subtle, and crucial for the note. If $X$ is projective  for example, when we say $X$ we mean the set of defining homogeneous polynomials in finitely many variables with coefficients in $\C$. This collection of polynomials is called  a {\it scheme}. On the other hand, only finitely many of those polynomials are necessary to describe them all (this is the Noetherian property of the ring of polynomials over a field), so in fact there is a ring $R$ of finite type over $\Z$ which contains all the coefficients. We write  $X_\C$ to remember $\C$, $X_R$ to remember $R$. 
We can then take any maximal ideal  $\frak{m}$ in $R$. The residue  field $R/\frak{m}$ is finite, say $\F_q$, and has characteristic $p>0$. Then we write $X_{\F_q}$ for the scheme defined by this collection of  polynomials  where the coefficients are taken modulo $\frak{m}$. Fixing an algebraic closure   $\F_q \subset \bar  \F_p$, and thinking of the polynomials as having coefficients in $\bar \F_q$ we write $X_{\bar \F_p}$ etc.

\medskip

When we say $X(\C)$, we mean the complex solutions of the defining polynomials.  (Of course there is the similar notion $X_R(R), X_{\F_q}( \F_q), X_{\bar \F_p}(\bar \F_p)=X_{\F_q}(\bar \F_p)$ etc.)

\medskip

 The notion of a quasi-projective complex variety $X$ is easily understood on its complex points $X(\C)$. They have to be of the shape  $\bar X(\C)\setminus Y(\C)$ where both $\bar X$ and $Y (\subset \bar X)$ are projective varieties.

 \medskip
 
We do not know how to characterize  the fundamental groups
$\pi_1(X(\C),x)$, where $x \in X(\C)$, among all possible $\pi_1(M, m)$.  In this small text, we 
use the following terminology:
\begin{defn}
A finitely presented group $\pi$ is said to come from geometry if it is isomorphic to $\pi_1(X(\C),x)$ where $X$ is a smooth connected quasi-projective  complex 
variety and $x\in X(\C)$.

\end{defn}
 {\it The aim of this note is to describe a few obstructions}  for a  finitely presented group to come from geometry. 
  
 \medskip
 
 {\it Acknowledgement:}  I thank  Johan de Jong, Michael Groechenig, and Moritz Kerz. 
  The material  in this  expository   note relies on joint  work or discussions  with them.  I thank Alexander Beilinson, Fran\c{c}ois Charles and Jakob Stix for friendly and helpful remarks on a first writing.  I thank the AMS for soliciting  this contribution.    I thank Han-Bom Moon for a careful reading and for suggestions aiming at making the note accessible to a wider audience.  The figures were drawn by him.
  
  \section{ Classical obstructions: topology and Hodge theory} \label{sec:Hodge}
 
\begin{figure}[!htb]
\begin{center}
\includegraphics[height=0.25\textheight]{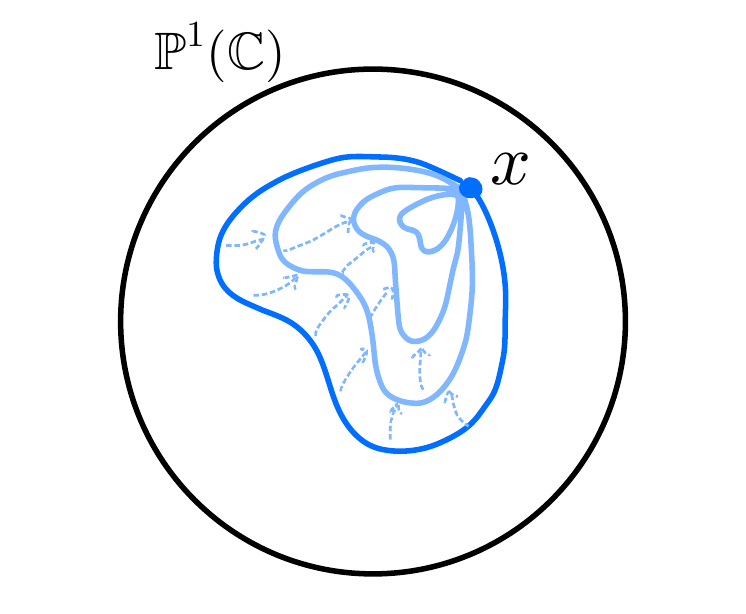}   
\end{center}
\caption{Any loop is retracted on $\P^{1}(\C)$}\label{fig:1}
\end{figure}

   A classical example comes from the uniformization theory of complex curves:  any {\it free} group on  $n$  letters, where $n$ is a natural number,   is the fundamental group of the complement of $(n+1)$-points on the Riemann sphere $\mathbb P^1$.     
     This is because  we understand exactly  $\pi_1(X(\C),x)$ if $X$ has dimension $1$, that is if $X(\C)$ is a Riemann surface. The simplest possible example  is  the Riemann sphere  $X=\P^1$. Then $\pi_1(X(\C),x)=\{1\}$ as any loop centered at $x$ can be retracted to a point, see Figure \ref{fig:1}. 

\newpage

\begin{figure}[!htb]
\begin{center}
\includegraphics[height=0.25\textheight]{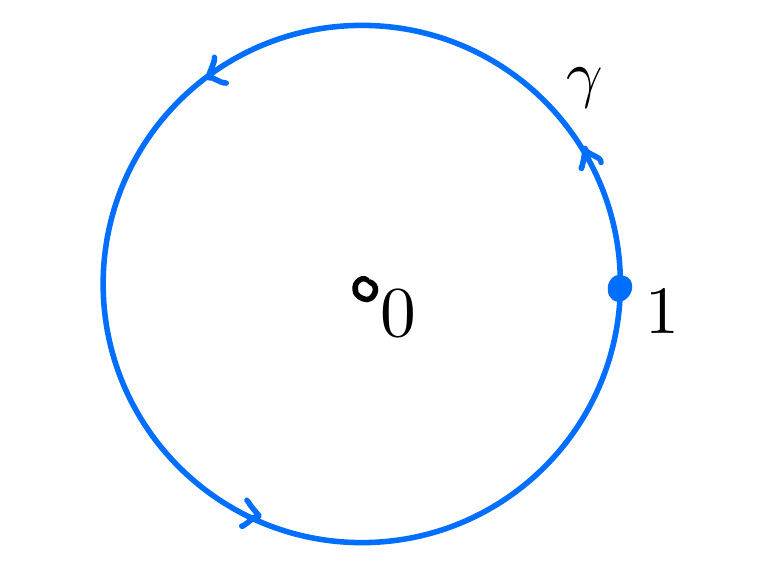}   
\end{center}
\caption{A nontrivial loop $\gamma$ on $\P^{1}(\C)\setminus\{0, \infty\}$}\label{fig:2}
\end{figure}

     The same holds true on $X=\A^1=\P^1\setminus \{\infty\}$.  The first interesting example is 
      $ X=\mathbb \P^1\setminus \{0,\infty\}$.  Then $X(\C)=\C\setminus \{0\}$, and $\pi_1(X(\C), 1)=\Z\cdot \gamma$ where $\gamma: [0 \ 1]\to \C\setminus \{0\}, \ t \mapsto {\rm exp}(2\pi \sqrt{-1}t)$ is the circle turning around the origin $\{0\}$, see Figure \ref{fig:2}. 
      
      More generally, 
     if a smooth compactification  $\bar X$ of $X$ has genus $g$, topologically $\bar X(\C)$ is a donut with $g$ holes. Then $\pi_1(\bar X(\C), x)$ is spanned by $2g$ elements $(a_i,b_i), i=1,\ldots, g$  with one relation $\prod_{i=1}^g [a_i, b_i]=1$. If $(\bar X \setminus X)(\C)$ consists of $(n+1)$ points, 
    $\pi_1(X(\C), x)$  is spanned by $2g+n+1$ elements $(a_i,b_i), i=1,\ldots, g, \ c_j, j=1,\ldots, n+1$ with one relation $\prod_{i=1}^g [a_i, b_i] \prod_{j=1}^{n+1} c_j=1.$ The literature is  full of beautiful coloured pictures visualizing this classical computation.


\medskip
   
 Beyond Riemann surfaces, that is,  for $X$ of dimension $\ge 2$, our understanding is very limited.

 \medskip

  The $2$ in the $2g$ in the previous example is more general:  by the fundamental structure theorem on finitely generated  $\Z$-modules,  the maximal abelian quotient $\pi_1(X(\C),x)^{\rm ab}$, that is, the abelianization of 
  $\pi_1(X(\C),x)$,
   is isomorphic to  a direct sum  of $ \Z^{\oplus b}$ for some natural number $b$ and of a finite abelian group $T$. 
   
   \medskip
   
\begin{figure}[!htb]
\begin{center}
\includegraphics[height=0.25\textheight]{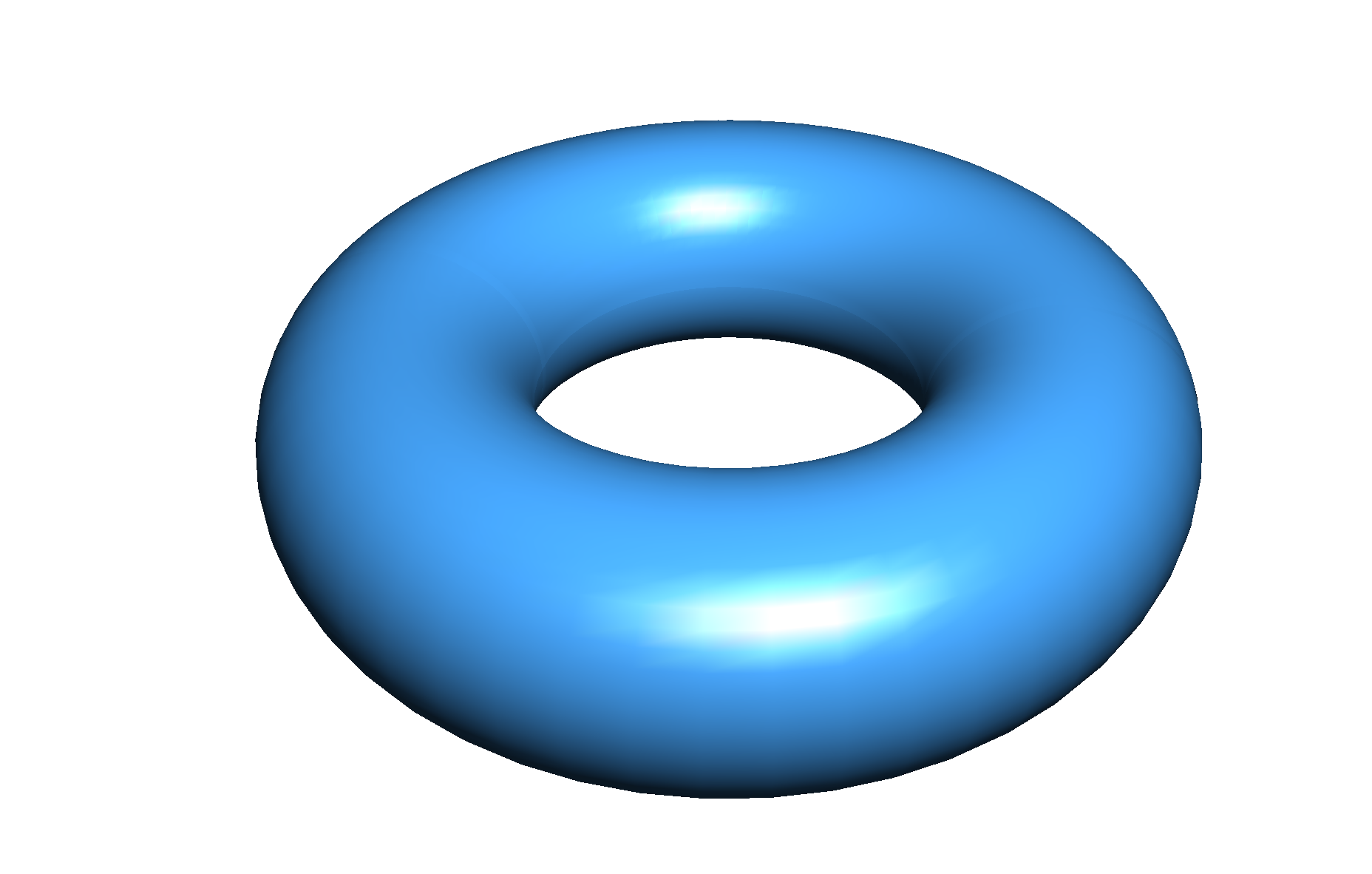}   
\end{center}
\caption{Riemann surface of genus $g = 1$}\label{fig:3}
\end{figure}

     {\it Any abelian finitely presented group $\Z^b\oplus T$ comes from geometry:}
   Serre's classical construction \cite{serre}  realizes any finite group as the fundamental group of the quotient $Z$ of a complete intersection of 
   large degree
   in the projective space of large dimension, while the fundamental group of  $
   (\P^1\setminus \{0,\infty\})(\C)$
    is equal to $\Z$, see Figure \ref{fig:2}. 
  As the  fundamental group of a product  variety is the product of   the fundamental groups of the factors (K\"unneth formula),  we can take   $X=  (\P^1\setminus \{0,\infty\})^b \times Z$  and there is {\it no} obstruction for $\Z^b\times T$ to be the 
   abelianization
   of the fundamental group of a  smooth connected quasi-projective complex variety. 
 If we require $X$ to be projective, then {\it Hodge theory},  more precisely, {\it  Hodge duality}  implies that $b$ is   {\it even}. This is the only obstruction as  we can then  take $X= E^{\frac{b}{2}} \times Z$ instead, where $E$ is any elliptic curve,   so $E(\C)$ is a donut with one hole,  so  $\pi_1(X(\C),x)=\Z^2$, see Figure \ref{fig:3}.

 \medskip
 
 In the same vein, but {\it much deeper} is the fact that the pro-nilpotent completion of $\pi_1(X(\C),x)$ (also called Mal\v{c}ev completion) is endowed with a {\it mixed Hodge structure}.  While so far we commented  the {\it topological} structure of $X(\C)$, Hodge theory studies in addition the {\it analysis} stemming   from  the {\it complex structure}, and the more refined properties, packaged in the notion of {\it K\"ahler geometry
 and harmonic theory}, which come from  the property that $X$ is defined  algebraically by complex polynomials.
 A modern way (due to Beilinson) to think of it  is to identify the Mal\v{c}ev completion with the cohomology of an (infinite) simplicial complex scheme and to apply the classical Hodge theory on its truncations. We do not elaborate further.

 \section{  Profinite completion: the \'etale fundamental group} \label{sec:et}

 Thus the difficulty lies in the kernel of the group to its abelianization. To study it, 
  we first introduce the   classical notion:
  
  \begin{defn} \label{defn:CLS}
 A   {\it complex local system} $\mathbb L_\rho$  is a complex linear representation $$\rho: \pi_1(X(\C),x)\to  \GL_r(\C),$$  considered {\it modulo  conjugacy} by $\GL_r(\C)$. The local system $\mathbb L_\rho$ is said to be irreducible if its underlying representation $\rho$ (thus defined modulo conjugacy) is irreducible. 
   \end{defn}
Why modulo conjugacy? A  path $\gamma_{xy}$ from $x$ to $y$ defines an isomorphism 
$\gamma_{yx}^{-1}\pi_1(X(\C),y)\gamma_{yx}=\pi_1(X(\C),x)$. This isomorphism is not unique, any other path from $x$ to $y$ differs from this one by left multiplication by a loop $\gamma_x \in \pi_1(X(\C),x)$ 
centered at $x$, which thus conjugates the isomorphism by $\gamma_x$.  Thus {\it not fixing the base point}
 forces us to consider representations modulo conjugacy.

\medskip

 As $\pi_1(X(\C),x)$ is finitely  presented, thus in particular finitely generated,   $\rho$ factors through 
 $\pi_1(X(\C),x)\xrightarrow{\rho_A} \GL_r(A)$ 
   where $A \subset \C$  is a ring of finite type. Any such $A$ 
   can be embedded into the ring of $\ell$-adic integers $\Z_\ell$ for some prime number $\ell$, 
  say $\iota: A \subset \Z_\ell$.  
   (For example if $A=\Z$,  $\iota$   has to be the natural pro-$\ell$-completion  for any choice of $\ell$. If $A=\Z[T]$ we take in $\Z_\ell$ a transcendental element over $\Q$ and send $T$ to it, etc. The main point is that the field of fractions $\Q_\ell$ of $\Z_\ell$, that is the field of $\ell$-adic numbers,   has infinite transcendence degree over $\Q$). 
  Thus the datum of $\rho$ is equivalent to the one of  $ \iota\circ \rho_A$ whose 
 range $\GL_r(\Z_\ell)$ is profinite. In particular,   $ \iota\circ \rho_A $ factors through the profinite completion  
 $${\frak{prof}}: \pi_1(X(\C),x)\to \pi_1(X(\C),x)^{\widehat{}} $$ and induces  $$\widehat{ \rho}: \pi_1(X(\C),x)^{\widehat{}} \to \GL_r(\Z_\ell),$$  a representation which is  {\it continuous} for the profinite topology on both sides. Recall that the profinite completion ${\frak{prof}}: \pi \to \pi^{\widehat{}} $
  of an abstract group $\pi$  is the projective limit over all finite quotients $\pi\to H$. It inherits the profinite topology  compatible with the group structure for which a basis  of open neighbourhoods of $1$ is defined to be the inverse image of $1 \in H$ by one of those projections.

  However, Toledo in  \cite{toledo}
 constructed a  smooth connected complex projective variety  $X$ with the property that 
  ${\frak{prof}}$
   is  {\it not} injective.   It answered a problem posed by  Serre. It  is an important fact which in particular  implies that  {\it the study of complex  local systems  ignores}  ${\rm Ker}(\frak{prof})$. 
  This leads us in two different directions.

 \medskip

   The invariants  $b$  and $T$ of the abelianization are seen on the complex abelian algebraic group 
 $${\rm Hom}(\pi_1(X(\C),x), \GL_1(\C))\cong (\C^\times)^b\times {\rm Hom}(T, \C^\times).$$
 Here   the notation $\C^\times$ means the set $\C\setminus \{0\}$ endowed with the (abelian) multiplicative group structure.
 
 \medskip
 
  More generally the finite generation of $\pi_1(X(\C),x)$  enables one to define a  ``moduli'' (parameter) space $M^{\rm irr} _B(X ,r)$ of all its irreducible local systems  $\mathbb L_\rho$ 
  in a given rank $r$. It is called the {\it Betti moduli space} of $X$  of irreducible local systems in rank $r$  or the {\it character variety} of $\pi_1(X(\C),x)$   of irreducible local systems in rank $r$. It is a complex   quasi-projective scheme of finite type.  Its study is the content of Simpson's {\it non-abelian Hodge theory} developed in \cite{simpson}. It is an {\it analytical  theory} relying on harmonic theory, as is classical Hodge theory.
  
 \medskip
  
  The second direction relies on the profinite completion homomorphism ${\frak{prof}}$.  
  By the Riemann existence theorem, a finite topological covering  is the complexification of a finite \'etale cover. Thus $\pi_1(X(\C),x)^{\widehat{}}$ is identified with the {\it \'etale fundamental group}  $\pi_1(X_\C,x)$ of the {\it scheme} $X_\C$ defined over $\C$, based at the complex point $x$, as defined by Grothendieck in \cite{SGA1}: 
  
  \medskip
  
   This profinite group is defined by its representations in finite sets. A representation  of $\pi_1(X_\C,x)$ in finite sets is ``the same'' (in the categorial sense) as  a pointed  (above $x$) finite \'etale cover of $X$.

   \medskip

  We denote   by 
  $$\rho_{\C, \ell}: \pi_1(X_\C,x)\xrightarrow{\widehat \rho} \GL_r(\Z_\ell) \hookrightarrow \GL_r(\bar \Q_\ell)$$ 
  the composite morphism.   Here $\bar \Q_\ell$ is an algebraic closure of $\Q_\ell$. 
  
  \medskip

  The notion of a complex local system (Definition~\ref{defn:CLS}) generalizes naturally: 
   \begin{defn} \label{defn:lLS}
 An  {\it $\ell$-adic  local system} $\mathbb L_{\rho_\ell} $ on the variety $X_\C$   is  a
continuous  linear representation $$\rho_\ell : \pi_1(X_\C,x)\to  \GL_r(\bar \Q_\ell),$$  considered {\it modulo  conjugacy} by $\GL_r( \bar \Q_\ell)$. The local system $\mathbb L_{\rho_\ell}$ is said to be irreducible if its underlying representation $\rho_\ell$ (thus defined modulo conjugacy) is irreducible. 
   \end{defn}

    As the kernel of the projection $\GL_r(\Z_\ell)\to \GL_r(\F_\ell)$ is a pro-$\ell$-group  (that is all its finite quotients $H$ have order of power of $\ell$), Grothendieck's  specialization's  theory  in {\it loc. cit.}   implies that 
    the specialization  homomorphism $$sp_{\C,\bar \F_p}: \pi_1(X_\C,x)\to \pi_1^t(X_{\bar \F_p},x)$$ induces an isomorphism on the image of $\rho_{\C, \ell}$ for $p$ larger than the order of  $\GL_r(\F_\ell)$. 
    Here
       $X_{\bar \F_p}$ is a reduction of $X_\C$  as explained in the introduction,  and is {\it good}, that is smooth, as well as the stratification  of the boundary divisor if $X$ is not projective. The upper script $t$ refers to the tame quotient of $\pi_1(X_{\bar \F_p},x)$ in case $X$ was not projective.  We do {\it not} detail  with precision the tameness concept, for which we refer to \cite{tame}.
    This roughly works as follows.  Representations in finite sets of the \'etale  fundamental group $ \pi_1(X_{\bar \F_p},x)$ which factor through the tame quotient  $\pi_1^t(X_{\bar \F_p},x)$
    have base change properties ``as if'' $X_{\bar \F_p}$  were proper. We can contract the fundamental group of $X$ over a $p$-adic ring $R$  with residue field $\bar \F_p$ to the one over $\bar \F_p$  in  the way we do topologically in order to identify the topological fundamental group of a tubular neighborhoud of a compact manifold  to the one of the compact manifold.       The natural identification of  $\pi_1(X_\C,x)$ with $\pi_1(X_K,x)$ where $K$ is an algebraic closure of the field of fractions of $R$ (this is called base change property) enables us  to define $sp_{\C,\bar \F_p}$.  Grothendieck computes that  $sp_{\C,\bar \F_p}$ induces an isomorphism on all finite quotients of 
   $  \pi_1(X_\C,x)$ and $\pi_1^t(X_{\bar \F_p},x)$  of order prime to $p$.

   \medskip

    The factorization defines the irreducible $\ell$-adic local system $\mathbb L_{\bar \F_p,\ell}$ on $X_{\bar \F_p}$ from which $\mathbb L_{\C,\ell}$ comes.       
      This leads us to study $\mathbb L_{\bar \F_p,\ell}$
  in order to derive {\it arithmetic properties} of  the initial $\mathbb L_\rho$. 
We can remark that again we know extremely little on the kernel of  $sp_{\C,\bar \F_p}$ and that the study of complex local systems ignores them as well, for a chosen $\iota: A\to \Z_\ell$ and $p$ large as before.

  \medskip
  
  On the other hand, $X_\C$ is defined over a field of finite type over $\Q$, thus with a huge Galois group, and $X_{\bar \F_p}$ is defined over a finite field $\F_q$ of characteristic $p>0$, with a very small Galois group isomorphic to $\hat \Z$, the profinite completion of $\Z$, topologically spanned by the Frobenius $\varphi$ of $\F_q$. Nonetheless, we shall see that this small Galois group yields non-trivial information. 
  
\medskip
 
 Our goal now is twofold. First
 we shall illustrate how to go back and forth between the Hodge theory side and the arithmetic side on a {\it particular example}.  This by far does not cover the whole deepness of the theory, but we hope that it gives some taste on how it functions.  Then we shall  mention on the way and at the end more general theorems to the effect that deep arithmetic properties stemming from the Langlands program,  notably the ``integrality''  illustrated on this particular example, enable one to find a new obstruction for the finitely presented group to come from geometry.

  \section{An example to study} \label{sec:ex}
  Let $X$ be a smooth   connected  quasi-projective  complex 
variety.  If $X$ is not projective, we fix a smooth projective compactification $X\hookrightarrow \bar X$ so that the divisor at infinity $D=\bar X\setminus X=\cup_{i=1}^M D_i $ is a strict normal crossings divisor (so its irreducible components $D_i$ are smooth and meet transversally).  For each $i$ we fix $r$ roots of unity $\mu_{ij}, j=1,\ldots, r$, possibly with multiplicity. They uniquely determine a conjugacy class  $T_i$ of a semi-simple matrix of finite order. 
The normal subgroup spanned by the conjugacy classes of small loops  $\gamma_i$ around the components $D_i$ is identified with  the kernel of the surjection $\pi_1(X(\C),x)\to \pi_1(\bar X(\C),x)$. 
We fix an extra natural number $\delta >0$.

\medskip

We make the following assumption \\ \ \\
{\bf Assumption} 
$(\star)_r$: \\ For a given rank $r \ge 2$, there are finitely many irreducible rank $r$ complex local systems $\mathbb L_\rho$  on $X$ such that the determinant of $\mathbb L_\rho$ has order dividing  $\delta$, and, if $X$ is not projective, such   that the semi-simplification of $\rho(\gamma_i)$ falls in $T_i$.

\medskip

 It is simple to describe $(\star)_r$:  in the Betti moduli space $M_B^{irr}(X,r)$ we have the subscheme  $M^{irr}_B(X,r, \delta, T_i)$ defined by the conditions $\{\delta, T_i\}$. The condition $(\star)_r$ means precisely that $M^{irr}_B(X,r, \delta, T_i)$  is $0$-dimensional, or equivalently that 
$M^{irr}_B(X,r, \delta, T_i)(\C)$  consists of finitely many points.

\medskip

Note  the condition on $\delta$  depends only on $\pi_1(X(\C),x)$ so could be expressed on the character variety, not however the condition on $T_i$. For this we have to know which $\gamma_i$ in $\pi_1(X(\C),x)$ come from the boundary divisor, so we need the geometry.

\medskip

If $r=1$, we drop  the condition on the determinant, and assume for simplicity that 
$X$ is projective. So the assumption becomes that there are finitely many irreducible rank $1$ complex local systems $\mathbb L_\rho$  on $X$.
This then  forces $b$ to be $0$, so $\pi_1(X(\C),x)^{\rm ab}$ to be finite.

\medskip

 Consequently, 
 those finitely many   $\mathbb L_\rho$ of rank $1$ have finite monodromy (i.e.  $\rho(\pi_1(X(\C),x)^{\rm ab})$ is finite).  
 This  implies that the $\mathbb L_\rho$  {\it come from geometry}, that is there is a smooth projective morphism $g: Y\to U\subset X$ where $U$ is a Zariski dense open in $X$ (in our case $U=X$), such that $\mathbb L_\rho$ restricted to $U$ is  a subquotient of the local system $R^ig_*\C$ coming from the representation of $\pi_1(U(\C),x)$ in $GL(H^i(g^{-1}(x),\C)))$ for some $i$ (in our case  $g$ is finite \'etale and $i=0$).
 
 \medskip
 
 A  different way of thinking of  finiteness is using 
 Kronecker's analytic criterion \cite{eilenberg}:  the set of the rank $1$ local systems is  invariant under the action of the automorphisms of $\C$ acting on $\GL_1(\C)=\C^\times$. 
 Finiteness of the monodromy is then equivalent to  the monodromy  being {\it unitary} (i.e. lying in $S^1\subset \GL_1(\mathbb C)=\mathbb C^\times$)
  and being {\it integral} 
  (i.e. lying  in $\GL_1(\bar \Z)\subset \GL_1(\C)$). 
We now discuss the generalization of these two  properties: {\it unitarity and integrality}. 

\medskip

 We first observe that $(\star)_r$  implies that the irreducible rank $r$ complex local systems are {\it rigid}  if we preserve the $\{\delta, T_i\}$ conditions.
As the terminology says, it means that we can not ``deform''  non-trivially the local system $\mathbb L_\rho$. Precisely it says that a formal deformation $$\rho_t: \pi_1(X(\C),x)\to \GL_r(\C[[t]])$$  of $\rho=\rho_{t=0}$  with the same $\{\delta, \mu_{ij}\}$ conditions does not move $\mathbb L_\rho$, that is there is a $g\in \GL_r( \C((t)))$    such that  in $\GL_r( \C((t)))$ the relation 
$$\rho_t=g\rho_{t=0} g^{-1}$$ 
holds.

\medskip

A classical example  where $(\star)_r$ is fulfilled is provided by Shimura varieties of real rank $\ge 2$.   Margulis super-rigidity  \cite{margulis} implies that all complex local systems are semi-simple and all irreducible ones are rigid.  {\it While  by super-rigidity they all are integral} (i.e. the image of the representations lie in $\GL_r(\bar \Z)$ up to conjugacy), {\it 
we do not know whether they come from geometry. }

\medskip

Another example is provided by connected smooth projective complex varieties $X$ with the property that all symmetric differential forms, except the functions,  are trivial. In this case, non-abelian Hodge theory 
implies 
$(\star)_r$ is fulfilled. Indeed, the Betti moduli space of semi-simple rank $r$ complex local systems is affine, while the moduli space of semi-stable Higgs bundles  with vanishing Chern classes (which we  discuss below) admits a projective morphism to the so-called Hitchin base. The latter   consists of one point under our assumption. As  by a deep theorem of Simpson~\cite{simpson}, both spaces are real analytically isomorphic, 
they are both affine and compact, thus are $0$-dimensional. 
It  is  proven in \cite{BKT}, using Hodge theory, the period domain and birational geometry, that  all  the $\mathbb L_\rho$ have then finite monodromy.  This yields a positive answer to a conjecture I had formulated.  As  the proof  uses Hodge theory, it is analytic. As of today, there is no arithmetic proof of the theorem.

\section{Non-abelian Hodge theory} \label{sec:simpson}
We first assume   that  $X$ is projective. We discuss a little more the notion of Higgs bundles mentioned above. 
Simpson in \cite{simpson}  constructs the moduli space  $M^{\rm s}_{Dol}(X,r, \delta)$ of stable Higgs bundles $(V, \theta)$  of degree $0$, where  $V$ is a vector bundle of rank $r$, $\theta: V\to \Omega^1_X\otimes V$ is a $\sO_X$-linear operator fulfilling the integrality  condition $\theta\wedge \theta=0$, such that ${\rm det}(V,\theta)$ has finite order dividing $\delta$.  (The integrality notion here is for the Higgs field $\theta$, and is not related to the integrality of a linear  representation mentioned in Section~\ref{sec:ex}). The stability condition is defined on 
the pairs $(V,\theta)$, that is one tests it on  Higgs subbundles. 
The finite order of ${\rm det}(V,\theta)$  implies that the underlying 
Higgs field of ${\rm det}(V,\theta)$ is 
equal  to $0$, so ${\rm det}(V,\theta) = ({\rm det}(V), 0).$ 
The moduli space $M^{\rm s}_{Dol}(X,r, \delta)$
is
 a complex scheme of finite type. 
It has  several features.

\medskip

 There is a real analytic isomorphism $M^{\rm irr} _{B}(X,r, \delta) \xrightarrow{\cong} M^{\rm s}_{Dol}(X,r, \delta)$.
So $(\star)_r$ implies   that $M^{\rm s}_{Dol}(X,r, \delta)$ consists of finitely many points.  

\medskip

Simpson defines on  Higgs bundles
the  algebraic $\C^\times$-action which assigns $(V, t\theta)$ to $(V, \theta)$ for $t\in \C^\times$. 
It preserves stability near $1\in \C^\times$ and semi-stability in general. Thus under the assumption $(\star)_r$, 
the $\C^\times$ -action stabilizes $M^{\rm s}_{Dol}(X,r, \delta)$ pointwise. 
Simpson proves in {\it loc. cit.} that  $\C^\times$-fixed points correspond  to polarized   complex  variations of Hodge structure (PCVHS). 
Mochizuki  in \cite{mochizuki} generalized this part of Simpson's theory to the smooth quasi-projective case so the conclusion remains valid in general. 

\medskip 
We summarize this section: 
{\it The assumption  $(\star)_r$ implies that the irreducible $\mathbb L_\rho$ of rank $r$, with determinant  of order diving $\delta$  and semi-simplification of $\rho(\gamma_i)$ falling  in $T_i$,  underlie a PCVHS. This property  is the analog of the unitary property in rank $r=1$}. 

\medskip

We can not expect more as already on Shimura varieties of  real rank $\ge 2$, not all local systems are unitary. If they all were, as they are integral, they would have finite monodromy. This is not the case.

\section{Arithmeticity} \label{sec:arithm}
 Again we fix $r$. Once we obtain the finitely many  local systems $\mathbb L_{\bar \F_p, \ell}$  on $X_{\bar \F_p}$ by specialization as in Section~\ref{sec:et}, also taking $p$  large enough so it be  prime to the orders of $\delta$ (and the $\mu_{ij}$  in case $X$ is not projective), we consider Grothendieck's homotopy exact sequence $$1\to \pi_1(X_{\bar \F_p},x)\to \pi_1(X_{\F_q},x)\to   \widehat{\Z}\cdot \varphi \to 1$$
 (\cite{SGA1}).
Here the finite field $\F_q\subset \bar \F_p$ is chosen so $X_{\bar \F_p}$ is defined over $\F_q$ and $\varphi$ is the Frobenius endomorphism of $\bar \F_p$ sending $\lambda$ to $\lambda^q$.

\medskip

Let us first discuss the meaning of the sequence in terms of finite \'etale covers. The surjectivity on the right  says that  if $ \F_q\subset \F_{q'}$ 
is a finite field extension, then the induced finite \'etale cover $X_{\F_{q'}}\to X_{\F_q}$ has no section.  The injectivity on the left says that any finite \'etale cover of $X_{\bar \F_p}$ can be dominating by one induced by  a finite \'etale cover of $X_{\F_q}$. The exactness in the middle says  that if  a finite \'etale cover of $X_{\F_q}$ acquires a section on $X_{\bar \F_p}$, then the induced cover of $X_{\bar \F_p}$  is completely split.

\medskip

The kernel of $\pi_1(X_{\bar \F_p},x)$ to its tame quotient  $\pi_1^t(X_{\bar \F_p},x)$ is normal in $\pi_1(X_{\F_q},x)$.  Thus a lift of  $\varphi$ to $\pi_1(X_{\F_q},x)$, which is well defined up to conjugation by $\pi_1(X_{\bar \F_p},x)$, 
acts by conjugation on $\pi_1^t(X_{\bar \F_p},x)$,  therefore  on   $\ell$-adic local systems on $X_{\bar \F_p}$  and respects tameness.

\medskip

 This action preserves  $r$,  irreducibility, $\delta$ and the $T_i$.  Thus $\varphi$ acts as a bijection on the finite set $\{\mathbb L_{\bar \F_p,\ell}\}$. 
  We conclude that replacing $q$ by some non-trivial finite power $q^{t}$, all $\mathbb L_{\bar \F_p,\ell}$ descend to $\ell$-adic  local systems 
$\mathbb L_{ \F_{q^t},\ell}$ 
on $X_{ \F_{q^t}}$. We say that the $\mathbb L_{\bar \F_p,\ell}$ are {\it arithmetic}. (This argument is adapted from \cite{integrality}). 

\medskip

We summarize this section: 
{\it The assumption  $(\star)_r$  implies that the local systems  $\mathbb L_{\bar \F_p,\ell}$ are arithmetic.} 

\medskip

 More generally, without the assumption $(\star)_r$ being fulfilled, 
Simpson proves in  \cite{simpson}  in all generality that the {\it  $\mathbb L_{\C,\ell}$ coming from irreducible rigid local systems are arithmetic}, that is they descend to $\ell$-adic local systems on $X_F$ where $F$  is a field of finite type over $\Q$.

\section{$\ell$-adic companions and integrality} \label{sec:comp}
Quoted  from \cite{eilenberg}, with adapted notation: 

\medskip

``Given a field automorphism $\tau$ of $\C$, we can postcompose the underlying monodromy representation of a complex local system $\mathbb L_\rho$ by  $\tau$ to define  a {\it conjugate}  complex local system  $\mathbb L_\rho^\tau$. Given a field automorphism $\sigma$ of $\bar \Q_\ell$, which then can only be continuous if it is the identity on $\mathbb \Q_\ell$, or more generally given a field isomorphism $\sigma$ between $\bar \Q_\ell$ and $\bar \Q_{\ell'}$  for some prime number $\ell'$, the postcomposition of a {\it continuous } non-finite monodromy representation is no longer continuous (unless $\ell=\ell'$ and $\sigma$ is the identity on $\Q_\ell$), so {\it we cannot define a conjugate} $\mathbb L^\sigma_{\C,\ell}$ of an $\ell$-adic local system by this simple postcomposition procedure.'' 

\medskip

However, when  $X_\C$  is replaced by  $X_{\bar \F_p}$, Deligne conjectured  in Weil II \cite{deligne} that we can. Let us  first state the conjecture.

\medskip

By the \v{C}ebotarev density theorem, an {\it irreducible} $\ell$-adic sheaf $\mathbb L_\ell$ defined by an irreducible continuous representation $\rho_\ell: \pi_1(X_{\F_q},x)\to \GL_r(\bar \Q_\ell)$ considered modulo conjugacy by $ \GL_r(\bar \Q_\ell)$,
 is determined {\it uniquely} by the characteristic polynomials $$P(\mathbb L_\rho, y,T)={\rm det} (T-\rho_\ell( Frob_y)) $$ for all closed point $y$ of $X_{\F_q}$, where $Frob_y$ is the arithmetic Frobenius at $y$.  This expression just means that the closed point $y$ has a residue field $\kappa(y) \subset \bar \F_p$ which is a finite extension of $\F_q$, say of degree $m_y$. Then  $y$ is a rational point of $X_{\kappa(y)}$, thus the
 conjugacy class of 
  $\varphi^{m_y}$ is well defined  as a subgroup  in $\pi_1(X_{\F_q},x)$.

 \medskip

The first part of the conjecture predicts that  if the determinant of $\rho$ has finite order, then
 $$P(\mathbb L_\ell, y,T)\in \bar \Q[T] \subset \bar \Q_\ell[T].$$
  In particular,   its $\sigma$-conjugate  $$P(y,T)^\sigma\in  \bar \Q[T] \subset \bar \Q_{\ell '} [T]$$ is defined.  
  
  \medskip
  
 The second part of the conjecture  predicts the existence of an irreducible $\ell'$-adic local system $\mathbb L^\sigma_{\ell}$ with  the property that 
$$P( \mathbb L^\sigma_{\ell }, y, T)= P(\mathbb L, y, T)^\sigma.$$  Again \v{C}ebotarev density theorem implies unicity up to isomorphism  once we know the existence.

\medskip

The two parts have been proven on smooth curves  $X_{\F_q}$  by Drinfeld in rank $r=2$ \cite{shtuka} and L. Lafforgue in any rank \cite{lafforgue} as a corollary of Langlands' conjecture over functions fields. It  thus uses  automorphic forms. Drinfeld's Shtukas also  imply that all 
$\mathbb L_\ell$ on $X_{\bar \F_q}$ come from geometry. 
 The existence of companions has been extended to smooth quasi-projective $X_{\bar \F_q}$  of any dimension by Drinfeld by arithmetic-geometric methods  \cite{drinfeld}, reducing to curves. The reduction method goes back initially to G\"otz Wiesend.
 
 \medskip
 
 Deligne's initial conjecture is for normal varieties. To go from smooth to normal varieties is still an open problem.

\medskip

So coming back to  $(\star)_r$,  we see that under this assumption, and abusing notation replacing $q^t$ by $q$, it holds that  for any $\ell'\neq \ell, \ell'\neq p$,  we have as many $\mathbb L^\sigma_{\F_q, {\ell}|_{X_{\bar \F_p}} }$ as $\mathbb L_{\bar \F_p, \ell}$. (The companion formation  preserves  the irreducibility on $X_{\bar \F_p}$,  the $T_i$ and $\delta$ as we had taken $p$ prime to all those orders, so in particular it  also preserves the  tameness).

\medskip

We now  use this fact to prove that {\it  the assumption $(\star)_r$ implies that all irreducible rank $r$ local systems with the conditions $ (\delta, \mu_{ij})$  are integral. }

\medskip 

The initial irreducible complex $\mathbb L_{\rho}$ being rigid, they are defined over $\sO_K[N^{-1}]$ where $K$ is a number field and $N$ is a positive natural number. 
 Precomposing  the $\mathbb L^\sigma_{\F_q, {\ell}|_{X_{\bar \F_p}} }$
  with $sp_{\C, \bar \F_p}$  and then by $\frak{ prof}$ yields 
as many  irreducible local systems as the $\mathbb L_{\rho}$ with the conditions given by $\delta$ and $\mu_{ij}$. They all are integral at the place above $\ell'$ determined by
 $$K\subset \bar \Q\subset \bar \Q_{\ell} \xrightarrow{\sigma} \bar \Q_{\ell '}.$$

\medskip

This construction is performed for all $\ell'\neq \ell$, $\ell'\neq p$,  all isomorphisms $\sigma$ between $\bar \Q_\ell$ and $\bar \Q_{\ell '}$, 
with $p$ large enough (larger than $N$, $\delta$, the order of the $\mu_{ij}$, and such that $sp_{\C, \bar \F_p}$ is well defined and surjective).   This finishes the proof of the integrality under the assumption $(\star)_r$. 

\medskip

This proof  is taken from \cite{integrality}.  It is shown in {\it loc. cit.} that the argument applies for cohomologically rigid local systems (a notion we do not detail here) without the assumption $(
\star)_r$. The assumption $(\star)_r$ does not imply cohomological rigidity, and cohomological rigidity does not imply $(\star)_r$. 

\section{ The $(\star)_r$ condition for an abstract finitely presented group} \label{sec:breuillard}
In \cite{dJE} we report   on the  example~\ref{ex:1} constructed by  Becker-Breuillard-Varlj\'u.
We quote from {\it loc. cit.} adapting the notation: 

\medskip

``
For $r=2$ and $\delta=1$, that is for ${\rm SL}_2$ representations, the authors compute that  $\Gamma_0$ has exactly  two  irreducible complex representations modulo conjugacy. 
The first one $\mathbb L_1$ is defined by
\ga{}{\rho_1(a)=\frac{1}{\sqrt{2}} \begin{pmatrix}
1 &  1 \\
-1 & 1 \notag
\end{pmatrix},
\ \ \rho_1(b)= \begin{pmatrix}
j &  0 \\
0 & j^2 
\end{pmatrix},
\notag} 
where $j$ is a primitive $3$-rd root of unity.  It is defined over $\Q(j)$. The local system  $\mathbb L_2$ is Galois conjugate to $\mathbb L_1$. The authors compute
\ga{}{ \rho(ab) =\frac{j}{\sqrt{2}} \begin{pmatrix}
1 &  j \\
-1 & j  \notag
\end{pmatrix}. }
As ${\rm Trace}(\rho(ab))=-\frac{1}{\sqrt{2}}$, $\mathbb L_1$ is not integral  at  $\ell=2$, so $\mathbb L_2$ is not integral at  $\ell=2$ either. Furthermore, 
$\rho_1(a)$ does not preserve the eigenvectors of $\rho_1(b)$, so $\mathbb L_1$ and thus $\mathbb L_2$ are irreducible with dense monodromy in ${\rm SL}_2(\C)$.
They also compute that those representations are  cohomologically rigid.''

\medskip 

So we see  that $\Gamma_0$ can not be  isomorphic to $\pi_1(X(\C),x)$ for a connected smooth {\it projective} complex variety $X$. We conclude that 
 {\it the  integrality  property in Section~\ref{sec:comp} is an  obstruction for a finitely presented group to come from projective geometry. }

\medskip

Jakob Stix remarks that $\Gamma_0^{\rm ab}$ is isomorphic to $\Z$, which has rank $1$, so $\Gamma_0$ obeys  the Hodge  theoretic obstruction mentioned in Section~\ref{sec:Hodge} as well. 

\medskip

The rest of the note is devoted to indicating how  to extend the obstruction based on integrality  to  {\it all}   connected quasi-projective 
varieties. 
 This is the content of \cite{dJE}. 

\section{de Jong's conjecture} \label{sec:dJconj}
If $X_{\bar \F_p}$ is a connected  normal   quasi-projective variety, and $\ell\neq p$  is a prime number,  de Jong conjectured in \cite{dejong} that an irreducible representation $$\pi_1(X_{\bar \F_p},x)\to \GL_r( \overline{\F_\ell((t))})$$ which is arithmetic 
 is in fact   constant in $t$, thus in particular has finite monodromy.  Here $\F_\ell((t))$ is the Laurent power series field over the finite field $\F_\ell$ and $\overline{\F_\ell((t))}$ is an algebraic closure.

 \medskip
 
  He shows that {\it assuming the conjecture,   irreducible representations $\pi_1(X_{\bar \F_p},x)\to \GL_r( \bar \F_\ell)$
always  lift to arithmetic $\ell$-adic  local systems if $X_{\bar \F_p}$ is a smooth connected curve}.  

\medskip

Drinfeld in  
\cite{kashiwara} applied this argument  to produce over  a connected normal complex quasi-projective variety $X_\C$ $\ell$-adic local systems with the property that via $sp_{\C,\bar \F_p}$ for $p$ large they are arithmetic over $X_{\bar \F_p}$. 

\medskip 

 de Jong's conjecture has been proved by B\"ockle-Khare in specific cases and Gaitsgory \cite{gaitsgory} in general  for $\ell \ge 3$. The latter proof  uses the geometric Langlands program.

\section{Weak integrality for groups}
Let $\Gamma$ be a finitely presented group, together with  natural numbers $r\ge 1, \delta\ge 1$. We define in \cite{dJE}  the following notion: $\Gamma$ has the {\it weak integrality property with respect to $(r,\delta)$} if, assuming there is an irreducible representation $\rho: \Gamma\to \GL_r(\C)$ 
with determinant of order $\delta$, then for any  prime number $\ell$, there is a representation $\rho_\ell: \Gamma\to \GL_r(\bar \Z_\ell)$  which is irreducible over $\bar \Q_\ell$ and of determinant of order  $\delta$. 

\medskip

The main theorem of {\it loc. cit.} is that {\it  if $X$ is a connected smooth quasi-projective complex variety, then $\Gamma=\pi_1(X(\C),x)$ does have the weak integrality property with respect to any $(r,\delta)$. }

\medskip

Using now the example by Becker-Breuillard-Varj\'u presented in  Section~\ref{sec:breuillard}, we see that their $\Gamma_0$ {\it does not  come from geometry} at all, whether the  desired  $X$ is assumed to be projective or 
quasi-projective.

\medskip

So we conclude that {\it  the weak integrality property for $\Gamma$  with respect to  all $(r,\delta)$ is an obstruction for $\Gamma$ to come from geometry. }  This new kind of obstruction does not rest on analytic methods, but on arithmetic properties, more specifically  { \it the arithmetic Langlands program for the existence of companions} and  {\it the geometric  Langlands program for de Jong's conjecture}, as we briefly discuss in the next and last section.

\section{Weak arithmeticity and density}

The main theorem of {\it loc. cit.} is proven by combining 
\begin{itemize}
\item[1)] the method discussed in Section~\ref{sec:comp} to show integrality once
 we have $\ell$-adic local systems on $X_{\bar \F_p}$ which are arithmetic;
 \item[2)]  and  the use de Jong's conjecture discussed in Section~\ref{sec:dJconj},  roughly as Drinfeld did  in \cite{kashiwara},   to produce many such arithmetic $\ell$-adic local systems on $X_{\bar \F_p}$.
 \end{itemize}
By Grothendieck's classical ``quasi-unipotent monodromy at infinity''  theorem \cite{serretate}, arithmetic tame $\ell$-adic local systems on $X_{\bar \F_p}$ 
have quasi-unipotent monodromies at infinity.  So  their pull-back to $X(\C)$ via $sp_{\C, \bar \F_p}$ and $\frak{prof}$ do as well.




\medskip 
In order  to apply the method described  in Section~\ref{sec:comp} involving the existence of $\ell$-companions ultimately  yielding  integrality, we need quasi-unipotent monodromies at infinity.  
 The method developed in \cite{dejong} shows that 
 those  in 
$M^{\rm irr}_B(X,r, {\rm torsion})$  are  Zariski dense,  where ``torsion'' refers to the determinant of $\mathbb L_\rho$ being torsion.   This is precisely this fact which enables one to ``forget'' the quasi-unipotent conditions at infinity and to develop the argument.

\

\bibliography{ExampleRefs}

\begin{bibdiv}
\begin{biblist}

\bib{BKT}{article}{
      author={Brunebarbe, Y.},
      author={Klingler, B.},
      author={Totaro, B.},
       title={Symmetric differential forms and the fundamental group},
        date={2013},
     journal={Duke Math. J.},
      volume={162},
       pages={2797\ndash 2813},
  url={https://projecteuclid.org/journals/duke-mathematical-journal/volume-162/issue-14/Symmetric-differentials-and-the-fundamental-group/10.1215/00127094-2381442.full},
      review={\MR{3127814}},
}

\bib{deligne}{article}{
      author={Deligne, P.},
       title={La conjecture de {W}eil {II}},
        date={1980},
     journal={Inst. Hautes \'Etudes Sci. Publ. Math.},
      volume={52},
       pages={137\ndash 252},
         url={http://www.numdam.org/item/?id=PMIHES_1980__52__137_0},
      review={\MR{0601520}},
}

\bib{dejong}{article}{
      author={de~Jong, J.},
       title={A conjecture on arithmetic fundamental groups},
        date={2001},
     journal={Israel J. Math.},
      volume={12},
       pages={61\ndash 84},
         url={https://link.springer.com/article/10.1007/BF02802496},
      review={\MR{1818381}},
}

\bib{dJE}{article}{
      author={de~Jong, J.},
      author={Esnault, H.},
       title={Integrality of the {B}etti moduli space},
        date={2023},
     journal={Trans. of the Am. Math. Soc.},
      volume={to appear},
       pages={18 pp},
         url={https://link.springer.com/article/10.1007/BF02802496},
}

\bib{kashiwara}{article}{
      author={Drinfeld, V.},
       title={On a conjecture of {K}ashiwara},
        date={2001},
     journal={Math. Res. L.},
      volume={8},
       pages={713\ndash 728},
  url={https://www.intlpress.com/site/pub/pages/journals/items/mrl/content/vols/0008/0006/a003},
      review={\MR{1879815}},
}

\bib{drinfeld}{article}{
      author={Drinfeld, V.},
       title={On a conjecture of {D}eligne},
        date={2012},
     journal={Mosc. Math. J.},
      volume={12},
       pages={515\ndash 542},
      review={\MR{3024821}},
}

\bib{shtuka}{article}{
      author={Drinfeld, V.},
       title={Langlands' conjecture for ${GL}(2)$ over function fields},
        date={1980},
     journal={Proceedings of the International Congress of Mathematicians
  (Helsinki, 1978)},
       pages={565\ndash 574},
      review={\MR{3024821}},
}

\bib{integrality}{article}{
      author={Esnault, H.},
      author={Groechenig, M.},
       title={Cohomologically rigid local systems and integrality},
        date={2018},
     journal={Selecta Math.},
      volume={24},
       pages={4279\ndash 4292},
         url={https://link.springer.com/article/10.1007/s00029-018-0409-z},
}

\bib{eilenberg}{article}{
      author={Esnault, H.},
       title={Lectures on {L}ocal {S}ystems in {A}lgebraic-{A}rithmetic
  {G}eometry},
        date={2023},
     journal={Lectures Notes in Mathematics, Springer Verlag},
      volume={2337},
       pages={70 pp},
}

\bib{gaitsgory}{article}{
      author={Gaitsgory, D.},
       title={On de {J}ong's conjecture},
        date={2007},
     journal={Israel J. Math.},
      volume={157},
       pages={155\ndash 191},
         url={https://link.springer.com/article/10.1007/s11856-006-0006-2},
      review={\MR{2342444}},
}

\bib{SGA1}{article}{
      author={Grothendieck, A.},
       title={Rev\^etements \'etales et groupe fondamental ({SGA}1)},
        date={1971},
     journal={Lectures Notes in Mathematics},
      volume={224},
       pages={xviii+327 pp},
      review={\MR{2017446}},
}

\bib{tame}{article}{
      author={Kerz, M.},
      author={Schmidt, A.},
       title={On different notions of tameness in arithmetic geometry},
        date={2010},
     journal={Math. Ann.},
      volume={346},
       pages={641\ndash 668},
         url={https://link.springer.com/article/10.1007/s00208-009-0409-6},
      review={\MR{2578565}},
}

\bib{lafforgue}{article}{
      author={Lafforgue, L.},
       title={Chtoucas de {D}rinfeld et correspondance de {L}anglands},
        date={2002},
     journal={Invent. Math.},
      volume={147},
       pages={1\ndash 241},
         url={https://link.springer.com/article/10.1007/s002220100174},
      review={\MR{1875184}},
}

\bib{margulis}{article}{
      author={Margulis, G.},
       title={Discrete subgroups of semisimple {L}ie groups},
        date={1991},
     journal={Ergebnisse der Mathematik und ihrer Grenzgebiete, Springer
  Verlag},
      volume={17},
       pages={x+388 pp},
         url={https://link.springer.com/article/10.1007/s002220100174},
      review={\MR{1090825}},
}

\bib{mochizuki}{article}{
      author={Mochizuki, T.},
       title={Kobayashi-{H}itchin correspondence for tame harmonic bundles and
  an application},
        date={2006},
     journal={Ast\'erisque},
      volume={309},
       pages={viii+117 pp},
      review={\MR{2310103}},
}

\bib{serre}{article}{
      author={Serre, J.-P.},
       title={Sur la topologie des vari\'et\'es en caract\'eristique $p$},
        date={1957},
     journal={Symposium internacional de topolog\'ia algebraica, Universidad
  Nacional Aut\'onoma de M\'exico and UNESCO, Mexico City},
       pages={24\ndash 53},
      review={\MR{0098097}},
}

\bib{simpson}{article}{
      author={Simpson, C.},
       title={Higgs bundles and local systems},
        date={1992},
     journal={Inst. Hautes \'Etudes Sci. Publ. Math.},
      volume={75},
       pages={5\ndash 95},
      review={\MR{1179076}},
}

\bib{serretate}{article}{
      author={Serre, J.-P.},
      author={Tate, J.},
       title={Good reduction of abelian varieties},
        date={1968},
     journal={Ann. of Math.},
      volume={88},
       pages={492\ndash 517},
         url={https://www.jstor.org/stable/1970722?origin=crossref},
      review={\MR{0236190}},
}

\bib{toledo}{article}{
      author={Toledo, D.},
       title={Projective varieties with non-residually finite fundamental
  groups},
        date={1993},
     journal={Inst. Hautes \'Etudes Sci. Publ. Math},
      volume={77},
       pages={103\ndash 119},
         url={http://www.numdam.org/item/?id=PMIHES_1993__77__103_0},
      review={\MR{1249171}},
}

\end{biblist}
\end{bibdiv}

\end{document}